\documentclass[11pt]{amsart}
\usepackage{amssymb}
\usepackage{graphicx}
\usepackage{tabularx}
\usepackage{epsfig,euler}
\usepackage{amsfonts}
\usepackage{amscd}
\usepackage{amsmath}
\newtheorem{defn}{Definition}[section]
\newtheorem{prop}{Proposition}[section]
\newtheorem{thm}[prop]{Theorem}
\newtheorem{cor}[prop]{Corollary}
\newtheorem{lemma}[prop]{Lemma}

\numberwithin{equation}{section}

\theoremstyle{remark}
\newtheorem{remark}{Remark}[section]

\newcommand{\Gr}{\operatorname{Gr}}

\newcommand{\C}{{\mathbb C}}
\newcommand{\Z}{{\mathbb Z}}

\begin{document}

\title[Sums of binomial determinants]{Sums of binomial determinants, non-intersecting lattice paths \\ and positivity of Chern-Schwartz-MacPherson classes}\author{Leonardo Constantin Mihalcea} \date{\today}
\begin{abstract}{We give a combinatorial interpretation of a certain positivity conjecture of Chern-Schwartz-MacPherson classes, as stated by P. Aluffi and the author in a previous paper. It translates into a positivity property for a sum of $p \times p$ determinants consisting of binomial coefficients, generalizing the classical Theorem of Lindstr\"om-Gessel-Viennot et al. which computes these determinants in terms of non-intersecting lattice paths. We prove this conjecture for $p=2,3$.}\end{abstract}
\maketitle
\section{Introduction} To any $2p$ lattice points $A_1, \cdots, A_p$ and $B_1, \cdots , B_p$ in $\Z^2$ one can associate a $p \times p$ matrix of binomial coefficients $M = (m_{ij})$, where $m_{ij} $ is equal to the number of lattice paths from $A_i $ to $B_j$, with each segment oriented either North-South, or West - East (see Figure 1 below). It is a classical result about binomial determinants (see e.g. \cite{L,GV2} or see \cite{K} and references therein) that if the points $A_i$ respectively $B_j$ are arranged, in order, from North-East to South-West, then the determinant of $M$ is non-negative, and counts $p-$tuples of non-intersecting lattice paths $\pi=(\pi_1, \cdots , \pi_p)$, where $\pi_i$ is a path from $A_i$ to $B_i$. In this note we conjecture a generalization of this result, where the conditions on the initial lattice points $A_i$ are relaxed (but $B_j$'s are still in the same NE-SW configuration). It is common for a fixed set of $A_i$'s in the relaxed hypothesis to yield a negative determinant, but when the {\em sum} of all the allowable configurations is considered, the result will be positive. The precise statement is given in Theorem \ref{thm:main} below. In this note we prove this conjecture for $p=2,3$.

Besides the intrinsic combinatorial interest, this conjecture has geometric significance: the sum of the determinants we consider is the coefficient of the fundamental class of a Schubert variety in the Grassmannian $\Gr(p,n)$ of $p-$planes in $\C^n$, for $n$ large enough, in the expansion of a Chern-Schwartz-MacPherson (CSM) class of a Schubert cell. In this geometric setting, the conjecture appeared in a previous paper by P. Aluffi and the author (\cite{AM}). For further details, we refer the reader to {\em loc. cit.}; the precise connection with the determinants considered herein is given in Remark \ref{rmk:csm} below. 

%In particular, a geometric proof of the positivity would amount to a proof of our conjecture.

\section{Statement of results} \subsection{Definitions and notations}\label{s:not} In this note a {\bf path} $\pi$ will be a lattice path in $\Z^2$ with the horizontal steps to the right and the vertical steps going down. For $A,B \in \Z^2$ the notation $\pi:A \to B$ means that $\pi$ starts at $A$ and ends at $B$. See the Figure 1 below.
 \begin{figure}[h!]
\begin{center}
\includegraphics[scale=0.5]{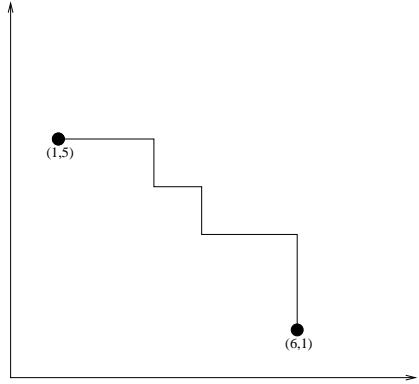}
\end{center}
\caption{A path from $(1,5)$ to $(6,1)$.}
%\label{figure_path1}
\end{figure}

By a {\bf partition} $\lambda$ we mean a decreasing sequence of nonnegative integers \[\lambda=(\lambda_1 \geqslant \lambda_2\geqslant \cdots \geqslant \lambda_p \geqslant 0).\]  Let $\lambda = (\lambda_1, \cdots , \lambda_p)$ and $\mu = (\mu_1, \cdots, \mu_p)$ be two such partitions. To this data we associate a family of $2p$ points $A_1(s), \cdots, A_p(s)$ and $B_1, \cdots B_p$ in $\Z^2$. The points $A_i(s)$ will depend on a sequence of parameters $s$ in a set $S$ which we define in the next paragraph. The points $B_j$, for $1 \leqslant j \leqslant p$ are defined simply by:
\[ B_j:=(\mu_j+p-j+1,\mu_j+p-j+1).\]

The set $S$ consists of sequences $s=(a_{i,j})$ of nonnegative integers indexed as the elements of a square matrix of order $p-1$, situated on or below the main diagonal:
 
 \begin{equation}\label{triangle} \begin{array}{cccc}
 a_{1,1} \\
 a_{2,1} & a_{2,2} \\
 \vdots & \vdots & \ddots \\
 a_{p-1,1}& a_{p-1,2} & \cdots & a_{p-1,p-1}
 \end{array} \end{equation}

\begin{defn}\label{deftriangular} Let $\lambda=(\lambda_1, \cdots , \lambda_p)$ be a partition. We say that the integers $a_{i,j}$ are in {\bf triangular order with respect to $\lambda$ } if:
\begin{enumerate}\item $0 \leqslant a_{i,j} \leqslant \lambda_{j+1}$. %are nonnegative.
\item The partial sums from the $j-th$ column, from row $j+1$ to $i+j$, for all $1 \leqslant i \leqslant p-1 - j $, are less than $a_{j,j}$, i.e.\begin{equation}\label{psum} a_{j+1,j}+...+a_{j+i,j} \leqslant a_{j,j}. \end{equation}
\end{enumerate}
\end{defn}
To simplify the notations in the upcoming formulae, we let $R_j$, for $2 \leqslant j \leqslant p-1$, respectively $C_j$, for $1 \leqslant j \leqslant p-2$ denote the partial sum on the $j-$th row, respectively $j-$th column of \ref{triangle}, excluding $a_{j,j}$:
\begin{equation}\label{rowsum} R_j:=a_{j,1}+...+a_{j,j-1}.\end{equation}
\begin{equation}\label{colsum} C_j:=a_{j+1,j}+...+a_{p-1,j}.\end{equation}
Set also $R_1=C_{p-1}=0$. Denote the set of all triangular sequences with respect to $\lambda$ by $S(\lambda)$ and an element of it by $s= (a_{ij})$.  We define the lattice points $A_j(s)$ by 
\[ A_j(s):=(p-j+1+a_{j,j}-R_j,\lambda_j+p-j+1-R_j)\] for $1 \leqslant j \leqslant p-1$.
If $j=p$, let the $x-$coordinate of $A_p(s)$ be \[x_{A_p(s)}:= 1+(C_1-a_{1,1})+(C_2-a_{2,2})+...+(C_{p-1}-a_{p-1,p-1})\] and the $y-$coordinate to be \[y_{A_p(s)}:=\lambda_p+x_{A_p(s)}.\] 

For $s \in S$ define the matrix $M(s) = (m_{ij}(s))$ by \[ m_{ij}(s) = \#\mathcal{P}(A_i(s) \to B_j) , \] where the right hand side denotes the number of paths from $A_i(s)$ to $B_j$. Let \[ c(\lambda,\mu) = \sum_{s \in S(\lambda)} \det M(s). \]

\subsection{Main result/conjecture} \begin{thm}[Positivity for $p \leqslant 3$ and general conjecture]\label{thm:main} The coefficient $c(\lambda,\mu)$ is positive if $p \leqslant 3$ and we conjecture it to be positive for all $p$. \end{thm}
\begin{remark} In the cases $p=2$ respectively $p=3$, positive combinatorial formulae for $c(\lambda,\mu)$ are given respectively in Corollaries \ref{posd2} and \ref{posd3} below. \end{remark}
\begin{remark}\label{rmk:csm} Choose $n$ large enough so that $\lambda_1, \mu_1 \leqslant n-p$. Denoting by $\mathbb{S}_\lambda^o$ respectively $\mathbb{S}_\mu$ the Schubert cell respectively the Schubert variety corresponding respectively to partitions $\lambda$ and $\mu$ (see \cite{AM} for details), $c(\lambda,\mu)$ is the coefficient of the CSM class of $\mathbb{S}_\lambda^o$ in the homology Schubert class $[\mathbb{S}_\mu]$. In the notation from \cite{AM},
$c(\lambda,\mu)=\gamma_{\underline{\lambda},\underline{\mu}}$. \end{remark}
\begin{remark}[Explicit definition of $M(s)$]  Given the triangular sequence $s=(a_{i,j})$, the matrix $M(s)$ is equal to:
 \begin{equation}\label{E:M(s)} M(s)= \begin{pmatrix}   \binom{\lambda_1-a_{1,1}}{\mu_1+R_1-a_{1,1}} & \binom{\lambda_1-a_{1,1}}{\mu_{2}-1+R_1 -a_{1,1}}&\cdots&\binom{\lambda_{1}-a_{1,1}}{\mu_p - (p-1)+R_1-a_{1,1}} \\    \binom{\lambda_2-a_{2,2}}{\mu_1+1+R_2-a_{2,2}} & \binom{\lambda_2-a_{2,2}}{\mu_{2}+R_2 -a_{2,2}}&\cdots&\binom{\lambda_{2}-a_{2,2}}{\mu_p - (p-2)+R_2-a_{2,2}} \\  \vdots & \vdots & \vdots & \vdots \\ \binom{\lambda_{p-1}-a_{p-1,p-1}}{\mu_1+p-2+ R_{p-1}-a_{p-1,p-1}} & 
\binom{\lambda_{p-1}-a_{p-1,p-1}}{\mu_{2}+p-3+R_{p-1}-a_{p-1,p-1}} &\cdots &   \binom{\lambda_{p-1}-a_{p-1,p-1}}{\mu_p-1+R_{p-1}-a_{p-1,p-1}} \\ \binom{\lambda_p}{\mu_1+p-1+\sum_{s=1}^{p-1} (a_{s,s}-C_s)} &\binom{\lambda_p}{\mu_2+p-2+\sum_{s=1}^{p-1} (a_{s,s}-C_s)} & \cdots &   \binom{\lambda_p}{\mu_p+ \sum_{s=1}^{p-1} (a_{s,s}-C_s)} \end{pmatrix}\end{equation}
i.e. the binomial coefficient on the row $r$ and column $c$, for $1 \leqslant r \leqslant p-1$ is equal to
\[\binom{\lambda_{r}-a_{r,r}}{\mu_{c}+r-c+R_{r}-a_{r,r}}.\]
For example, in the case $p=2$, the triangular sequence $S(\lambda)$ consists of all $(a_{11})$ such that $0 \leqslant a_{11} \leqslant \lambda_2$ and $M(a_{11})$ is given by:
%\[ \gamma_{\underline{\lambda},\underline{\mu}} =
% \sum_{i_1^1=\lambda_2} \det \begin{pmatrix} \binom{\lambda_2}{\gamma_2+k} &  \binom{\lambda_2}{\gamma_1+1+k} \\ 
%    \binom{\lambda_1-k}{\gamma_2-1-k} & \binom{\lambda_1-k}{\gamma_1-k}   \end{pmatrix}
 
 \begin{equation}\label{E:det_p=2} \begin{pmatrix} \binom{\lambda_1-a_{1,1}}{\mu_1-a_{1,1}} &  \binom{\lambda_1-a_{1,1}}{\mu_2-1-a_{1,1}} \\ \binom{\lambda_2}{\mu_1+1+a_{1,1}} & \binom{\lambda_2}{\mu_2+a_{1,1}}\end{pmatrix}.\end{equation}
Similarly, in the case $p=3$, the triangular sequences consist of triples $(a_{21},a_{22},a_{11})$\begin{footnote}{We used the ordering $(a_{21},a_{22},a_{11})$ rather than $(a_{11},a_{12},a_{22})$ to be consistent with the notation $(a_{21},a_{22},a_{11}) = (i,j,k)$ used throughout the paper starting from the next paragraph.}\end{footnote} such that
 \[ 0 \leqslant a_{21} \leqslant a_{11}; \textrm{ } 0 \leqslant a_{22} \leqslant \lambda_3; \textrm{ } 0 \leqslant a_{11} \leqslant \lambda_2 , \] and $M(s)$ is:
%\begin{equation}c(\lambda;\gamma)=\sum_{k=0}^{\lambda_2}\sum_{i=0}^k\sum_{j=0}^{\lambda_3} \det   \begin{pmatrix}
%    \binom{\lambda_3}{\gamma_3+(k-i)+j} &  \binom{\lambda_3}{\gamma_2+1+(k-i)+j} &  \binom{\lambda_3}{\gamma_1+2+(k-i)+j} \\ 
%    \binom{\lambda_2-j}{\gamma_3-1+i-j} & \binom{\lambda_2-j}{\gamma_2+i-j} & \binom{\lambda_2-j}{\gamma_1+1+i-j} \\ 
%    \binom{\lambda_1-k}{\gamma_3-2-k} &     \binom{\lambda_1-k}{\gamma_2-1-k}  &     \binom{\lambda_1-k}{\gamma_1-k}  
%  \end{pmatrix}\end{equation}
\begin{equation}\label{sumd3} \begin{pmatrix} 
 \binom{\lambda_1-a_{1,1}}{\mu_1-a_{1,1}} &     \binom{\lambda_1-a_{1,1}}{\mu_2-1-a_{1,1}}  &     \binom{\lambda_1-a_{1,1}}{\mu_3-2-a_{1,1}} \\   \binom{\lambda_2-a_{2,2}}{\mu_1+1+a_{2,1}-a_{2,2}} & \binom{\lambda_2-a_{2,2}}{\mu_2+a_{2,1}-a_{2,2}} & \binom{\lambda_2-a_{2,2}}{\mu_3-1+a_{2,1}-a_{2,2}} \\  \binom{\lambda_3}{\mu_1+2+(a_{1,1}-a_{2,1})+a_{2,2}} &  \binom{\lambda_3}{\mu_2+1+(a_{1,1}-a_{2,1})+a_{2,2}} &  \binom{\lambda_3}{\mu_3+(a_{1,1}-a_{2,1})+a_{2,2}}\end{pmatrix}.\end{equation}
\end{remark}
\begin{remark} The positivity conjecture for $c(\lambda,\mu)$ was checked on the computer for all pairs $\lambda,\mu$ included in the partitions $7^5 = (7,7,7,7,7), 5^6=(5,5,5,5,5,5),10^4=(10,10,10,10)$ etc.  \end{remark}

\subsection{An example for $p=3$.}
Let $\lambda=(3,3,3)$ and $\mu=(2,2,1)$. To avoid carrying subscripts in the case $p=3$, we identify the triangular sequence $(a_{21},a_{22},a_{11})$ to $(i,j,k)$, so that 
 \[ 0 \leqslant k \leqslant \lambda_2; \textrm{ } 0 \leqslant i \leqslant k; \textrm{ } 0 \leqslant j \leqslant \lambda_3 \quad . \] The lattice points $A_\ell,B_\ell$, for $A_\ell = A_\ell(i,j,k)$, $1 \leqslant \ell \leqslant 3$ are given by
\[A_1 = (k+3,\lambda_1+3);\quad A_2=(2-i+j,\lambda_2+2-i);\quad A_3:=(1-k+i-j,\lambda_3+1-k+i-j), \] \[B_1 = (5,5), B_2=(4,4),B_3=(2,2) \quad . \]
%and the coefficient $c(\lambda,\mu)$ has the formula:
%\begin{equation}\label{det_intro} c(\lambda,\mu)=\sum_{k=0}^{\lambda_2}\sum_{i=0}^k\sum_{j=0}^{\lambda_3} \det M(i,j,k)\end{equation} where $M(i,j,k)$ is the matrix
%\begin{equation}\label{intro_sumd3}  %\gamma_{\underline{\lambda},\underline{\mu}} =
% %\sum_{k=0}^{\lambda_2}\sum_{i=0}^{k}\sum_{j=0}^{\lambda_3} \det 
% \begin{pmatrix} 
% \binom{\lambda_1-k}{\mu_1-k} &     \binom{\lambda_1-k}{\mu_2-1-k}  &     \binom{\lambda_1-k}{\mu_3-2-k} \\   \binom{\lambda_2-j}{\mu_1+1+i-j} & \binom{\lambda_2-j}{\mu_2+i-j} & \binom{\lambda_2-j}{\mu_3-1+i-j} \\  \binom{\lambda_3}{\mu_1+2+(k-i)+j} &  \binom{\lambda_3}{\mu_2+1+(k-i)+j} &  \binom{\lambda_3}{\mu_3+(k-i)+j}\end{pmatrix}.\end{equation}
Using the the version of Lindstr\"om-Gessel-Viennot from \cite{GV2}, Thm. 1, it follows that each of the determinants of matrices $M(i,j,k)$ counts {\em signed} triples of non-intersecting lattice paths. We are forced to include signs since a non-intersecting triple may also arise from a permutation $w \in Sym(3)$ of the initial points. In this case a triple $\pi_w=(\pi_1,\pi_2,\pi_3)$ where $\pi_\ell:A_{w(\ell)} \to B_\ell$ has to be counted with the sign $\varepsilon(w) = (-1)^{l(w)}$, where $l(w)$ is the length of $w$. The content of the Theorem \ref{thm:main} is that all triples counted negatively are cancelled by the positive ones.   

In fact, we will prove more: if the sum $j+k$ is fixed, say $j+k=f$ then \begin{equation}\label{E:partsum} c(\lambda,\mu;f) := \sum_{(i,j,k) \in S(\lambda), j+k=f} M(i,j,k) \end{equation} is non-negative and there exists an $f$ such that this sum is positive. As an example, let $j+k=2$ (so $f=2$). The configurations arising from this situation are those from Figure 2 below. Then $c((3,3,3),(2,2,1);2)$ is the sum of $6$ determinants, and it can be written as: \[ c((3,3,3),(2,2,1);2) = 0 + 3 + 6 + 3 + 3 - 3 = 12 \quad .\] 
\begin{figure}[h!]\label{configjk}
\includegraphics[scale=.3]{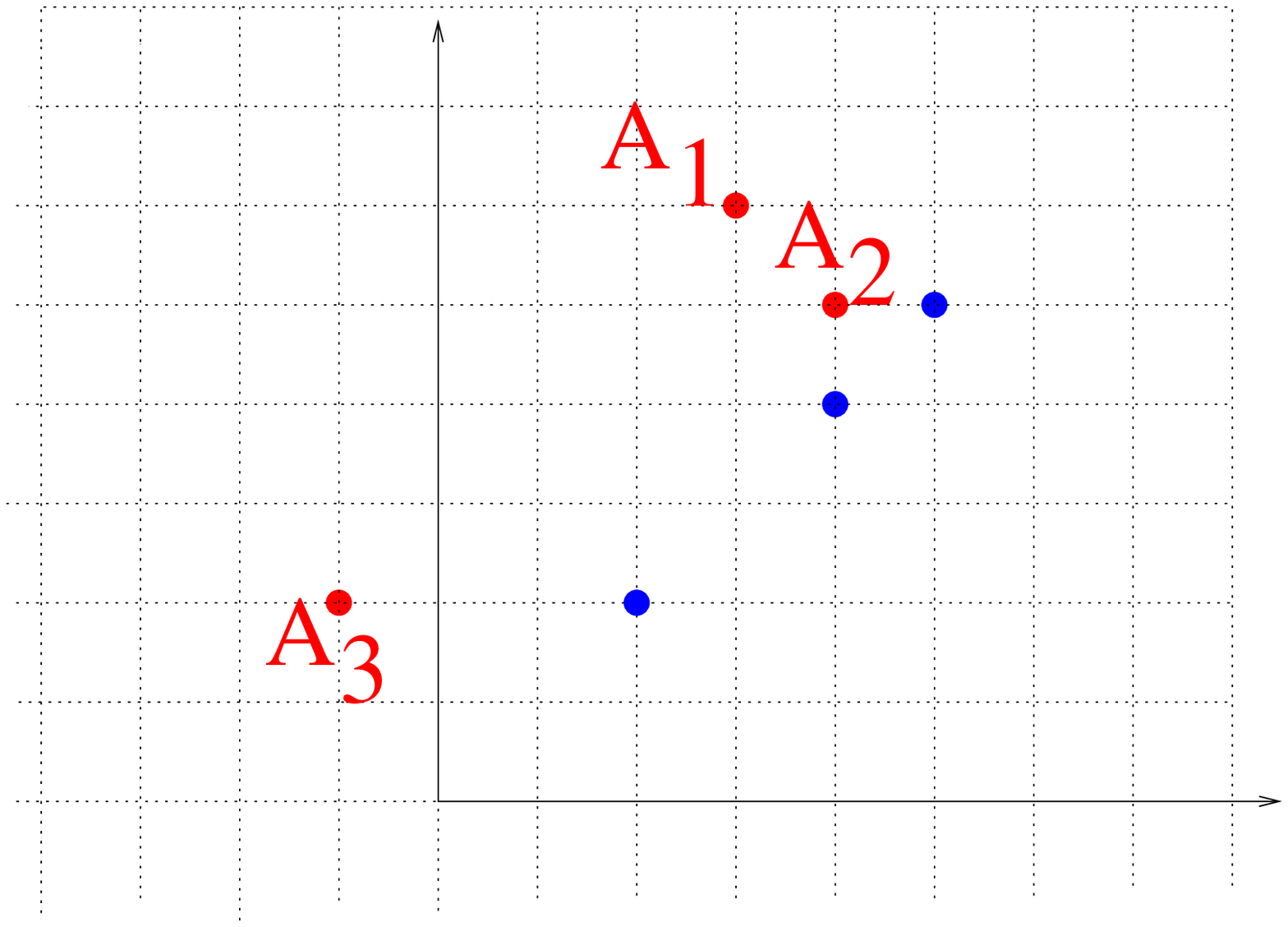}
\includegraphics[scale=.3]{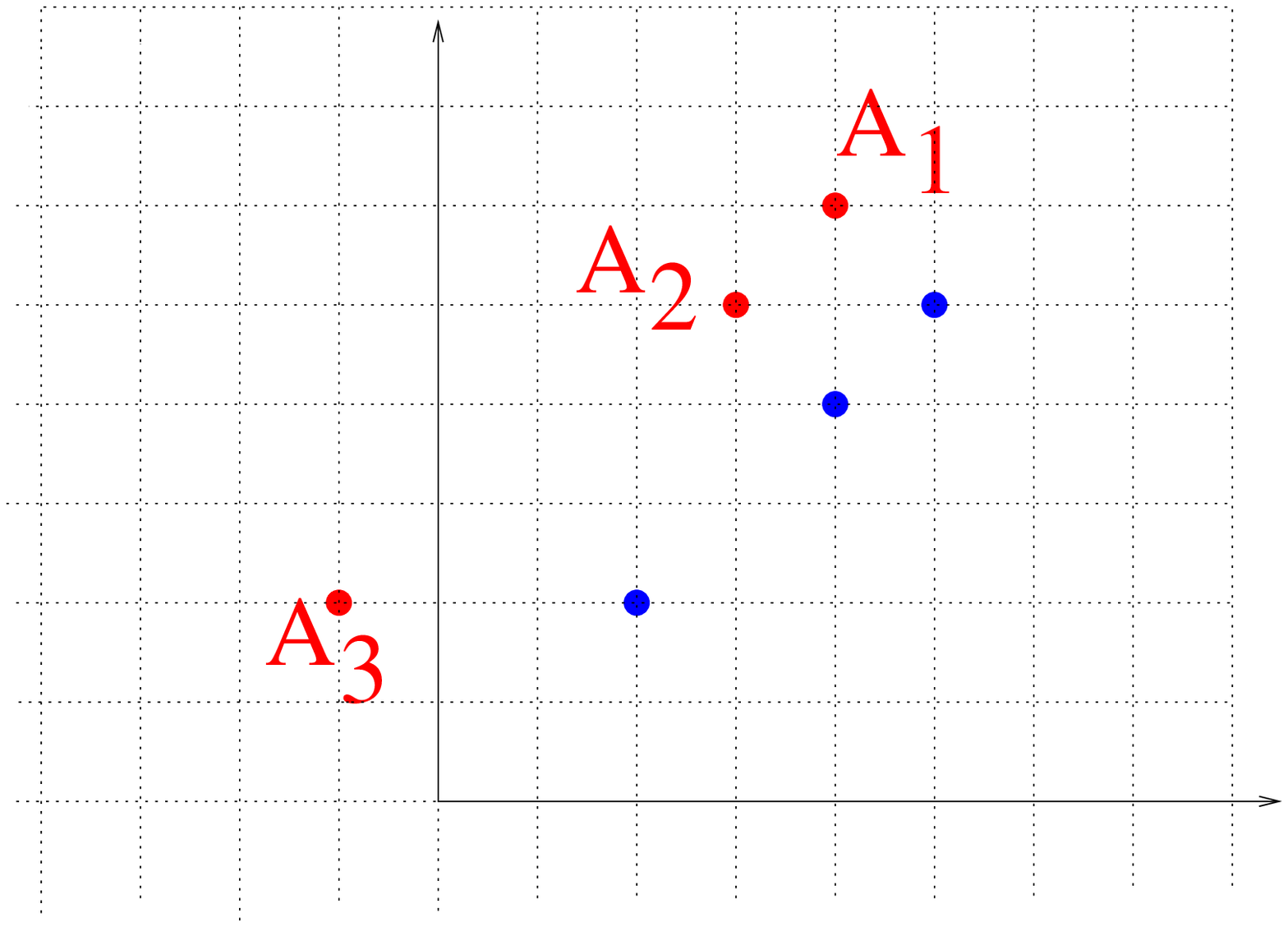}
\includegraphics[scale=.3]{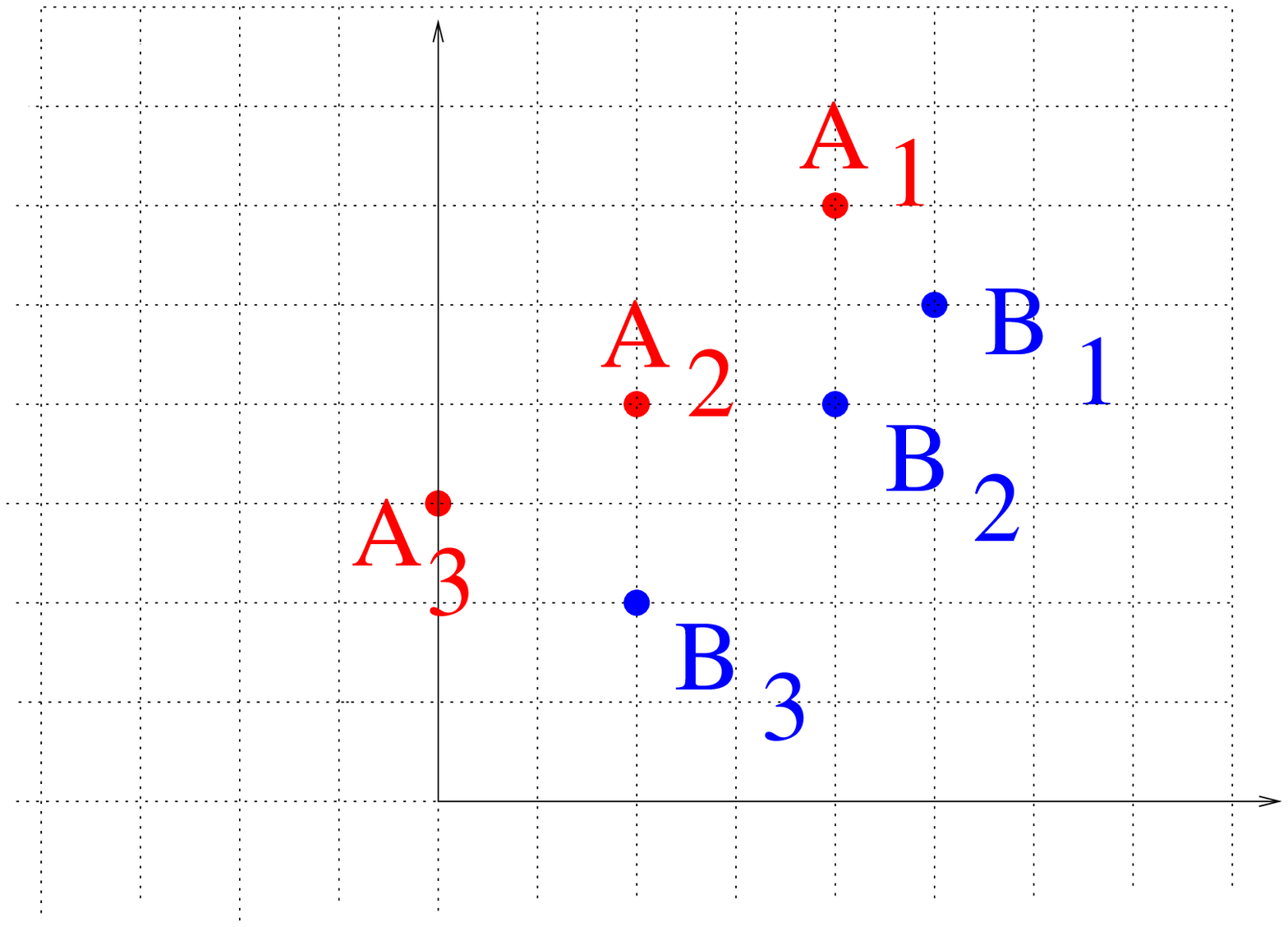}
\includegraphics[scale=.3]{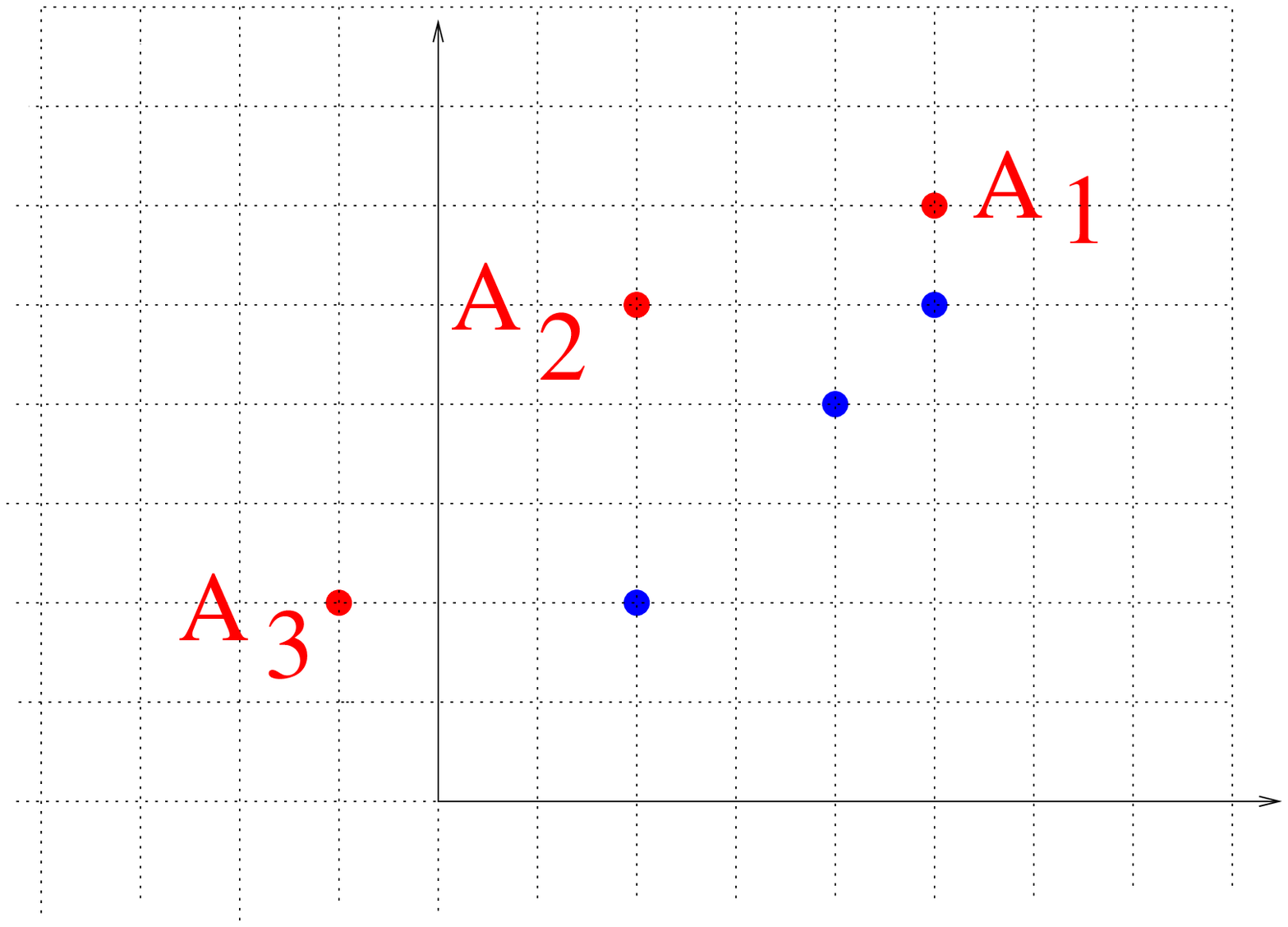}
\includegraphics[scale=.3]{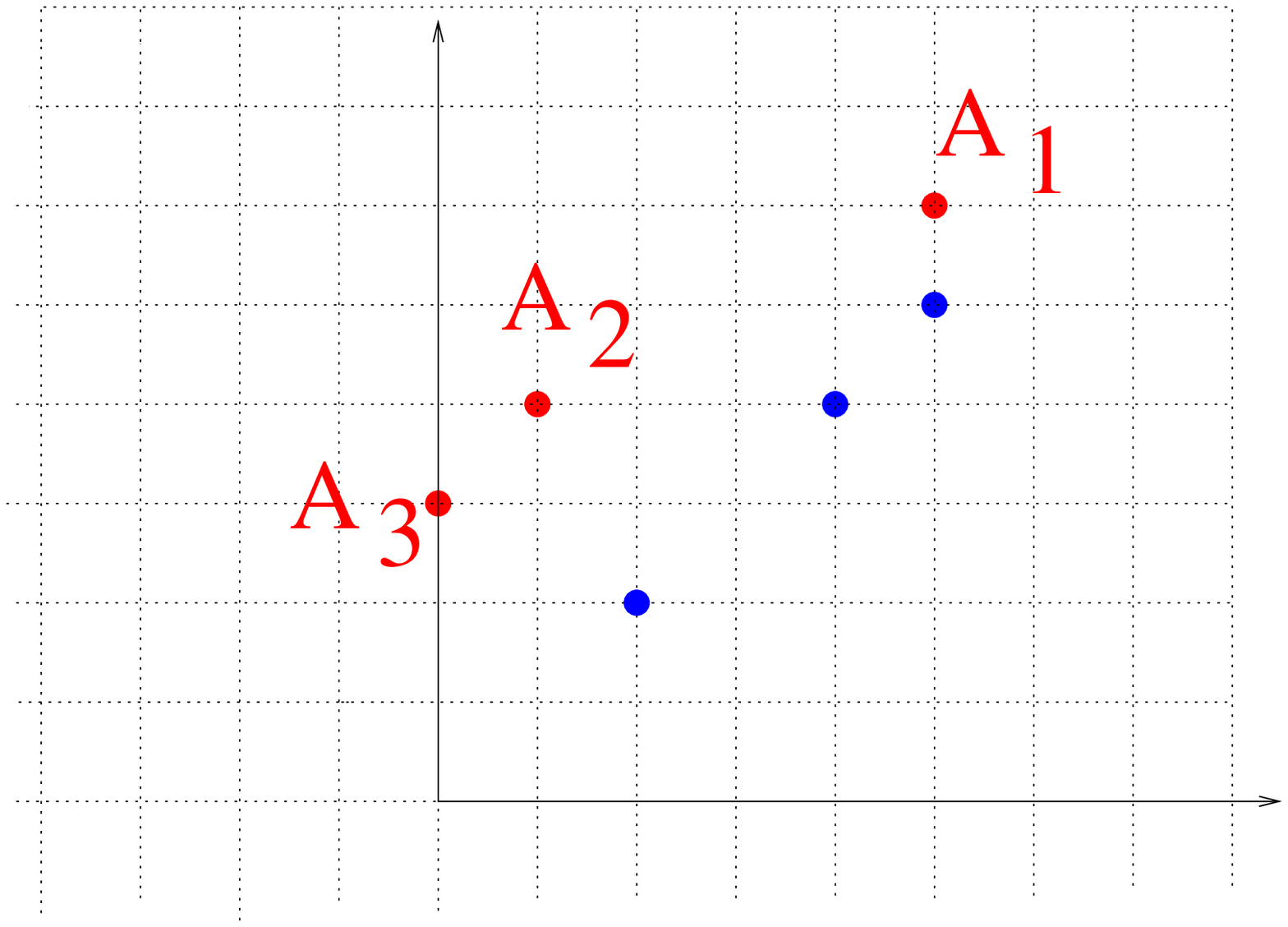}
\includegraphics[scale=.3]{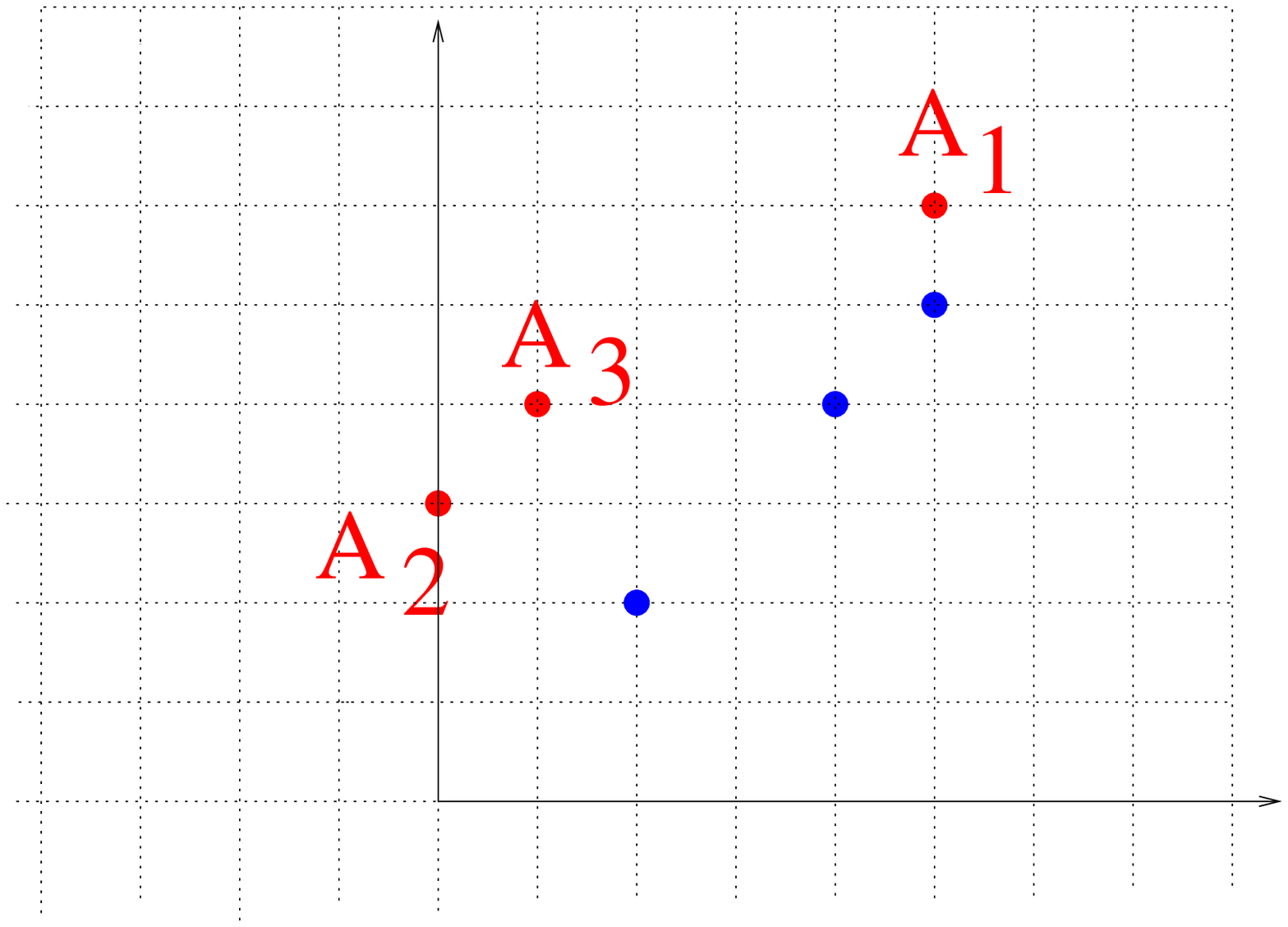}
\caption{Configurations of A's and B's for $\lambda=(3,3,3),\mu=(2,2,1)$ and $j+k=2$; figures correspond, left-right, top-down to $(0,2,0),(0,1,1),(1,1,1),(0,0,2),(1,0,2),(2,0,2)$. The blue dots represent the points $B_\ell$, which do not vary with $(i,j,k)$.}
\end{figure}

Note, for example, that $\det M(0,2,0)=1-1=0$ since there is one triple of non-intersecting paths,  $(A_1,A_2,A_3) \to (B_1,B_2,B_3)$ counted with $+1$, and one triple $(A_2,A_1,A_3) \to (B_1,B_2,B_3)$ counted negatively. This is different from the case $(i,j,k) = (2,0,2)$ when all the triples are counted negatively (in fact, this is just an ordinary Gessel-Viennot determinant, with the second and the third row swapped). 
 
\subsection{Idea of proof} To show that the coefficient $c((\lambda_1,\lambda_2,\lambda_3),(\mu_1,\mu_2,\mu_3);f)$ is nonnegative it is enough to prove that there is an injective map from the non-intersecting triples which count negatively to those counting positively. We show how this map is constructed for the negative paths computing $\det M(0,2,0)$ and $\det M(2,0,2)$ in Figure 2 above. To shorten notations, we denote by $\mathcal{P}((A_1,A_2,A_3) \longrightarrow (B_1,B_2,B_3))$ the set of {\em non-intersecting} triples of lattice paths $\pi=(\pi_1, \pi_2, \pi_3)$ where $\pi_\ell:A_\ell \to B_\ell$.

The map will distinguish between two cases: one when $w$ is the transposition $(12)$ and
one when $w=(23)$; $(i,j,k)=(0,2,0)$ corresponds to the first case, while $(2,0,2)$ to the second. We will show among other things, in \S \ref{s:pos3}, that these are the only configurations resulting in (non-intersecting) triples counted negatively.

If $(i,j,k)=(0,2,0)$, from a triple of paths \[(\pi_1,\pi_2,\pi_3) \in \mathcal{P}((A_2,A_1,A_3) \longrightarrow (B_1,B_2,B_3)),\] we construct a triple \[(\pi_1^*,\pi_2^*,\pi_3^*) \in \mathcal{P}((A_1^*,A_2^*,A_3^*) \longrightarrow (B_1,B_2,B_3))\] where $(A_1^*,A_2^*,A_3^*) $ is the triple corresponding to $(i^*,j^*,k^*)=(0,1,1)$; the path $\pi_3$ remains unchanged, so $\pi_3^*=\pi_3$. As for $\pi_1^*$ respectively $\pi_2^*$, they are constructed using certain `surgery' on $\pi_1$ and $\pi_2$ respectively. This process, described below, is shown in Figure 3. First, one translates the source $A_2$ of $\pi_1$ horizontally to left, say $x$ units, until it hits $\pi_2$. Let $A_2^*$ be this intersection point and define $\pi_2^*$ to be the portion of $\pi_2$ starting at $A_2^*$.  Similarly, given the $x$ units from the previous step, one translates the portion of $\pi_2$ from $A_1$ to $A_2^*$ horizontally to the right $x$ units, and form the new path $\pi_1^*$. %This process is illustrated in figure \ref{inv12}. 
\begin{figure}[h!]\label{inv12}
\includegraphics[scale=.5]{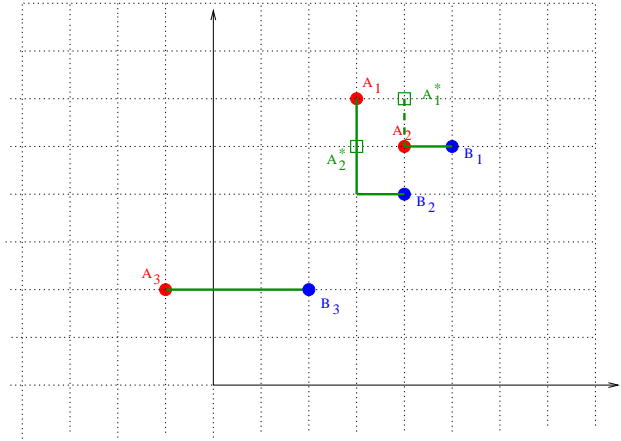}
\caption{Construction of paths $\pi_1^*$  and $\pi_2^*$ corresponding to the inversion $(12)$}\end{figure}
Note that in this case, the triple $(i^*,j^*,k^*)$ corresponding to $(A_1^*,A_2^*,A_3^*)$ is obtained from the initial $(i,j,k)$ by making \begin{equation}\label{tr1}i^*:=i,\textrm{   } j^*:=j-x,\textrm{   } k^*:=k+x, \end{equation} and such a transformation leaves $S(\lambda)$ and the sum $j+k$ invariant, provided that $x$ is small enough. 

A similar procedure, using now a {\em diagonal translation} with slope $1$, can be used to construct a positive triple out of one corresponding to the inversion $(23)$. This is illustrated in Figure 4. \begin{figure}[h!]\label{inv23} 
\includegraphics[scale=.5]{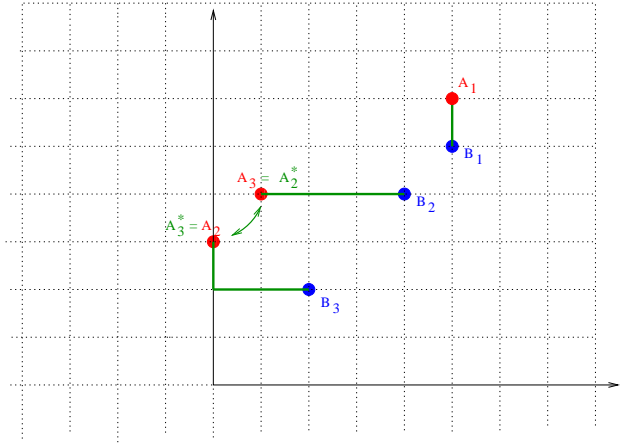}
\caption{Construction of paths $\pi_2^*$  and $\pi_3^*$ corresponding to the inversion $(23)$}\end{figure} In this case, the newly obtained triple
$ (A_1^*,A_2^*,A_3^*)$ is via the transformation \begin{equation}\label{tr2}i^*:=i-x,\textrm{   } j^*:=j,\textrm{   } k^*:=k \quad , \end{equation} which again preserves $S(\lambda)$ and the sum $j+k$ for small $x$.\\
%It is a coincidence that, in this case, $a_2^*=a_3$ and $a_3^*=a_2$. A more general configuration can be seen on the figure ? below. 

{\em Acknowledgments.} I am grateful to Paolo Aluffi for some inspiring conversations and to Christian Krattenthaler for his valuable suggestions and pointing out useful references in the area. Most of this work was done while enjoying the hospitality and support of Max-Planck Institute f\"ur Mathematik, Bonn.
\section{Proof of the main result}
\subsection{Preliminaries on non-intersecting lattice paths} We use the notations from \S \ref{s:not}. We fix two partitions $\lambda=(\lambda_1, \cdots , \lambda_p)$ and $\mu = (\mu_1, \cdots , \mu_p)$. The following is the connection between the geometric determinantal formulae for 
CSM classes from Theorem 3.4 in \cite{AM} and the determinants considered in this paper. Recall that the determinant $M(s)$ was defined in equation (\ref{E:M(s)}).
\begin{prop}\label{p:connection} Let $s\in S(\lambda)$ be a triangular sequence and let $A_\ell(s), B_\ell$ be the associated lattice points. Then \[ M(s) = \bigl( \# \mathcal{P}(A_r(s) \to B_c)\bigr)_{1 \leqslant r,c \leqslant p} \quad . \] \end{prop}
\begin{proof} This is a straightforward computation, taking into account that the number of paths between the lattice points $A=(a_1,a_2)$ and $B=(b,b)$, such that $A$ is NW of $B$ is $\binom{a_2 - a_1}{b-a_1}$.\end{proof}   

We recall next the (unweighted) version of the classical theorem of Lindstrom-Gessel-Viennot which allows any configuration of the initial points and end points. Recall that $\mathcal{P}(E \to F)$ denotes the set of lattice paths from $E$ to $F$.
 
\begin{thm}[Thm. 1 in \cite{GV2}]\label{thmnonint} Let $E_i,F_j$ be $2p$ lattice points, with $1 \leqslant i,j \leqslant p$. Then the determinant $\det ( \# \mathcal{P}(E_i \to F_j)_{ 1 \leqslant i,j \leqslant p}$ is equal to $\sum_{\pi_w} \varepsilon(w)$, where $w$ is a permutation in $Sym(p)$, $\varepsilon(w)$ is its signature and the sum is over all $p-$tuples of paths  \[\pi_w=(\pi_1^w,...,\pi_p^w)\] with $\pi_i^w:E_{w(i)} \to F_i $ , such that no two paths $\pi_i^w$ and $\pi_j^w$ intersect.\end{thm}
%\begin{proof} See \end{proof}

\subsection{Possible configurations for the points $A_\ell(s)$ and $B_\ell$.} Since the parts of the partition $\mu$ are decreasing it follows that for $\ell_1>\ell_2$, the point $B_{\ell_1}$ is strictly North-East of $B_{\ell_2}$ (see figure below).
\begin{figure}[h!]
\begin{center}
\includegraphics[scale=0.5]{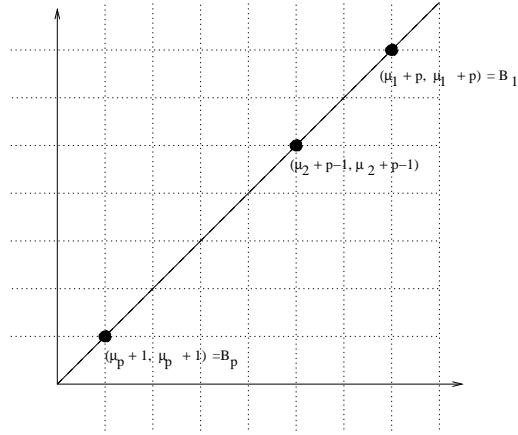}
\end{center}
\caption{The end points $B_\ell$.}
%\label{figure_path1}
\end{figure}
The points $A_\ell=A_\ell(s)$, for a fixed triangular sequence $s$, are not arranged as nicely. However, the following holds:
\begin{lemma}\label{lemmaconfig}(a) Let $s$ be a triangular sequence. Then $A_1(s)$ is strictly North and strictly East of $A_p(s)$, i.e. $x_{A_1(s)}>x_{A_p(s)}$ and $y_{A_1(s)}>y_{A_p(s)}$. \\
(b) $A_1(s)$ is strictly North of $A_\ell(s)$, for all $2 \leqslant \ell \leqslant p-1$.
\end{lemma}
\begin{proof} This is a straightforward computation.\end{proof}

\subsection{Two distance functions between paths} Let $\pi$ and $\pi'$ be respectively two paths between $A$ and $B$, and $A'$ and $B'$. We define two distances between $\pi$ and $\pi'$. One is a horizontal distance, denoted $D_h(\pi,\pi')$ and the other is a diagonal distance, denoted $D_d(\pi,\pi')$. These distances will be the quantities $x$ used to define the transformations described in equations (\ref{tr1}) and (\ref{tr2}).

We define $D_h$ first. This distance will only be defined provided that $y_{A}>y_{A'}$, i.e. that $A'$ is strictly to the South of $A$. The reader can refer to Figure 3, with $(A,B)=(A_1,B_2)$ and $(A',B')=(A_2,B_1)$. Assume that a horizontal line $L'$ passing through $A'$ intersects $\pi$ in a point $C$. Then $D_h(\pi,\pi')=l$ where $l$ is the length of the segment $A'C$. If $L'$ doesn't intersect $\pi$ then we set $D_h(\pi,\pi')=\infty$. In Figure 3, $D_h=1$.

To define $D_d$, we take diagonal lines $L$ and $L'$ of slope $1$ starting respectively from $A$ and $A'$ in both directions (in Figure 4, $(A,B)=(A_3,B_2), (A',B')=(A_2,B_3)$, and the lines $L$ and $L'$ happen to coincide). If at least one of $L$ or $L'$ intersects respectively $\pi'$ and $\pi$ (say in $C'$ or $C$), then $D_d(\pi,\pi')$ is the (diagonal) length of the segment $AC'$ or $A'C$ (necessarily just one) formed in this way. If neither of $L$ and $L'$ intersects the associated paths, define $D_d(\pi,\pi')=\infty$. The following is immediate:

\begin{lemma}\label{lemmadistance} Assume that $\pi_r:A_r \to B_r$, $1 \leqslant r \leqslant 2$ are two paths, and that $B_1$ and $B_2$ are on the main diagonal. Then $D_d(\pi,\pi')<\infty$.\end{lemma}

\subsection{Two swaps}
Let $\pi_1:A_1 \to B_2$ and $\pi_2:A_2 \to B_1$ be two paths as in Figure 3. We will define two `swaps' between $\pi_1$ and $\pi_2$; one horizontal and one diagonal provided that the corresponding distance between the paths is not infinity. The result will be a new pair of paths, $(\pi_1^*,\pi_2^*)$ and new sets of points $A_r^*$, $r=1,2$. The end points $B_r$ remain fixed. 

In fact, we only define the {\em horizontal swap} as in Figure 3 above and we let the reader to fill in the details for the diagonal swap, using the Figure 4. Let $A_r=(x_{A_r},y_{A_r})$, $r=1,2$.

Assume that $D_h(\pi_1,\pi_2)=l$ and that this distance is realized by the segment $A_2A_2^*$, with $A_2^* \in \pi_1$. Then let $\pi_2^*$ to be the partial path obtained from $\pi_1$ by chopping off the part from $A_1$ to $A_2^*$. To define $\pi_1^*$ we translate horizontally the path from $A_1$ to $A_2^*$ and attach it to $\pi_2$  such that $A_2^*$ becomes $A_2$. Then
$A_1^*$ will be the new starting point of $\pi_1^*$. In terms of coordinates, if $(x_r,y_r)$ and $(x_r^*,y_r^*)$ are the coordinates respectively of $A_r$ and $A_r^*$, then
\[ (x_1^*,y_1^*)=(x_1\pm l,y_1) \textrm{ and } (x_2^*,y_2^*)=(x_2\mp l,y_2) \] where $\pm$ is decided by the orientation of the segment $A_2A_2^*$: plus if $x_{A_2^*} < x_{A_2}$, minus otherwise. Similarly, if $D_d(\pi_1,\pi_2)=l'$ then $ (x_r^*,y_r^*)=(x_r\pm l',y_r \pm l') $.
\subsection{Positivity for $p=2$}
Recall that in the case $p=2$ a triangular sequence $s$ has just one element, denote it $i$, satisfying $0 \leqslant i \leqslant \lambda_2$. Then the situation is very simple: for any $i$, the lattice points $A_1(i),A_2(i),B_1,B_2$ satisfy the hypothesis of the Gessel-Viennot Theorem, i.e. $A_1(i)$ is North-East of $A_2(i)$ for all $i$. Therefore the intermediate determinants from equation (\ref{E:det_p=2}) are non-negative. To state the precise formula, let $\Pi(i)$ denote the set of all {\em non-intersecting} pairs of paths $(\pi_1,\pi_2)$, with $\pi_r: A_r(i) \to B_r$.
\begin{cor}[Positivity for $p= 2$]\label{posd2} Let $\lambda=(\lambda_1,\lambda_2)$ be a partition. Then for any partition $\mu$, the coefficient $c(\lambda,\mu)$ is equal to \[ \sum_{i=0}^{\lambda_2}\#\Pi(i) \quad .\]\end{cor} 
%\begin{proof} Fix $i$ such that $0 \leqslant i \leqslant \lambda_2$ and let $a_1,a_2$ be its associated lattice points. Since $a_1$ is strictly NE of the point $a_2$ (i.e. $x_{a_1}>x_{a_2}$ and $y_{a_1}>y_{a_2}$), by lemma \ref{lemmaconfig}, any path from $a_1$ to $B_2$ must intersect any path from $a_2$ to $B_1$. Then the condition in the hypothesis of the Gessel-Viennot Theorem is satisfied, therefore the determinant of $M((i))$ is non-negative. To finish the proof, it is enough to observe that the zero-sequence, i.e. $i=0$, clearly triangularly ordered with respect to $\lambda$, corresponds to a determinant which is positive, since there is at least one path from $a_r$ to $B_r$ in this instance.\end{proof}

\subsection{Positivity for $p=3$}\label{s:pos3} To shorten notations, let, as before, $(a_{21},a_{22},a_{11})= (i,j,k)$.
%\[k:=a_{11}; \qquad j:=a_{22}; \qquad i:=a_{21}.\]
The constraints for the triangular sequence with respect to $\lambda$ translate to:
\begin{equation}\label{*}0 \leqslant k \leqslant \lambda_2; \qquad  0 \leqslant i \leqslant k; \qquad 0 \leqslant j \leqslant \lambda_3.\end{equation}
Fix such a sequence $(i,j,k)$; recall that \[A_1 = (k+3,\lambda_1+3);\quad A_2=(2-i+j,\lambda_2+2-i);\quad A_3:=(1-k+i-j,\lambda_3+1-k+i-j).\] Fix also a partition $\mu=(\mu_1 \geqslant \mu_2 \geqslant \mu_3)$ which determines $B_1,B_2,B_3$:
\[ B_r:=(\mu_r+3-r+1,\mu_r+3-r+1), \quad 1 \leqslant r \leqslant 3 .\]

Invoking again the classical Gessel-Viennot theorem, if $A_1,A_2,A_3$ are arranged, in order, (weakly) from NE to SW, the corresponding determinant will be nonnegative. Since $A_1$ is always strictly NE of $A_3$, and strictly North of $A_2$ (by Lemma \ref{lemmaconfig}) it follows that there are only two possibilities which may yield a negative
determinant:\\

{\bf Case 1.} {\em $A_2$ is strictly S and strictly E of $A_1$ and there is a triple of nonintersecting paths $\pi=(\pi_1,\pi_2,\pi_3)$ such that $\pi_1:A_1 \to B_2$, $\pi_2:A_2 \to B_1$ and $\pi_3:A_3 \to B_3$} (see e.g. the configuration corresponding to $(i,j,k)=(0,2,0)$ in Figure 2).
In this case, let $l:=D_h(\pi_1,\pi_2)$ be the horizontal distance between $\pi_1$ and $\pi_2$. Clearly $l < \infty$. Then perform the horizontal swap to $(\pi_1,\pi_2)$ to define $(\pi_1^*,\pi_2^*)$, with $\pi_r^*:A_r^* \to B_r$. The new starting points $A_1^*$ and $A_2^*$ are given by the sequence $(i^*,j^*,k^*)$ where\begin{equation}\label{eqcoordc1} i^* = i; \qquad j^*=j-l; \qquad k^*=k+l.\end{equation}  Note that $\pi_3$ is not affected; denote by $\pi^*=(\pi_1^*,\pi_2^*,\pi_3)$ the new obtained triple. We have to show that $(i^*,j^*,k^*)$ satisfy the constraints conditions in \ref{*}. Indeed, the relative position of the paths $\pi_1$ and $\pi_2$ implies that $x_{A_2}-l \geqslant x_{A_1}$, i.e. that \[ 2+i-j-l\geqslant k+3. \] In particular \[ j-l \geqslant k+1-i >0 \textrm{ and } k+l \leqslant i-j-1 < \lambda_2.\] Another important feature of the swaps performed is that the paths in the new triple $\pi^*$ are nonintersecting. This follows immediately from their  construction: $A_3$ is SW of $A_1$ and $\pi_1^*$ is obtained by moving the `head' of $\pi_1$, containing $A_1$, horizontally to the right.

{\bf Case 2.} {\em $A_2$ is {\bf not} weakly North {\bf and}  weakly East of $A_3$ and there is a triple of nonintersecting paths $\pi=(\pi_1,\pi_2,\pi_3)$ such that $\pi_1:A_1 \to B_1$, $\pi_2:A_2 \to B_3$ and $\pi_3:A_3 \to B_2$}. There are three situations, according to $A_2$ being NW, SE or SW of $A_3$ (see subcases 2.1, 2.2 and 2.3 below). In all three situations one performs a diagonal swap to $(\pi_2,\pi_3)$. We obtain a new triple $\pi^*=(\pi_1,\pi_2^*,\pi_3^*)$, with $\pi_r^*:A_r^* \to B_r$ and the new sequence $(i^*,j^*,k^*)$ defining $A_r^*$ ($r=2,3$) is given by: \begin{equation}\label{eqcoordc2} i^*= i-l; \qquad j^*=j; \qquad k^*=k.\end{equation}  As in Case 1, we have to show that $(i^*,j^*,k^*)$ satisfies the constraints $\ref{*}$, i.e. that $i \geqslant l$.

{\em Subcase 2.1. $A_2$ is NW of $A_3$, i.e. $x_{A_2} < x_{A_3}$ and $y_{A_2}\geqslant y_{A_3}$.}  In this case \[l \leqslant x_{A_3}-x_{A_2}=(1-k+i-j)-(2-i+j)=2i-k-j-1\] which shows that
\[ i - l \geqslant k-i+j+1 >0.\]

{\em Subcase 2.2. $A_2$ is SE of $A_3$, i.e. $x_{A_2} \geqslant x_{A_3}$ and $y_{A_2}< y_{A_3}$.}  In this case \[l \leqslant y_{A_3}-y_{A_2}=(1+\lambda_3-k+i-j)-(2+\lambda_2-i)=(\lambda_3-\lambda_2)+2i-k-j-1\] which shows that
\[ i - l \geqslant (\lambda_2-\lambda_3)+k-i+j+1 >0.\]

{\em Subcase 2.3. $A_2$ is SW of $A_3$, i.e. $x_{A_2} \leqslant x_{A_3}$ and $y_{A_2}\leqslant  y_{A_3}$, but $A_2 \neq A_3$.}  In this case \[l \leqslant \max \{x_{A_3}-x_{A_2},y_{A_3}-y_{A_2}\}\] and the computation reduces to one from Subcases 2.1 or 2.2 above. 

We also need to show that the diagonal swap produces a non-intersecting triple of paths. This is done separately for each of the subcases above, and it should be clear from the construction.\\

To finally show positivity, let $c(\lambda,\mu;f)$ be the partial sum from equation (\ref{E:partsum}) defining the coefficients of CSM classes in the case $p=3$, obtained by fixing $k+j=f$.
\begin{thm}\label{T:posd3} Let $0 \leqslant f \leqslant \lambda_2+\lambda_3$. Then the partial sum $c(\lambda,\mu;f)$ is nonnegative and $c(\lambda,\mu;0)>0$ if $\mu \subset \lambda$ (i.e. $\mu_r \leqslant \lambda_r$, for $ 1 \leqslant r \leqslant 3$). \end{thm}
{\em Proof.} The second part of the Theorem is immediate: if $f=0$, then $(i,j,k)=(0,0,0)$, and the points $A_1,A_2,A_3$ are arranged from NW to SE, thus satisfying the hypothesis of the Gessel-Viennot Theorem. Since $\mu \subset \lambda$, each $B_\ell$ is SE of $A_\ell$, so there is at least one non-intersecting triple of paths $\pi:(A_1,A_2,A_3) \to (B_1,B_2,B_3)$. To prove the first part it is enough to show that there cannot be triples $\pi=(\pi_1,\pi_2,\pi_3)$ and $\pi^{**}=(\pi_1^{**},\pi_2^{**},\pi_3^{**})$ such that:
\begin{enumerate} \item $\pi$ and $\pi^{**}$ are triples of nonintersecting paths.

\item $\pi$ creates a $(12)$ inversion, as in Case 1 above, i.e. $\pi_1:A_1(i,j,k) \to B_2$, $\pi_2:A_2(i,j,k) \to B_1$, $\pi_3:A_3(i,j,k) \to B_3$.
\item  $\pi^{**}$ creates a $(23)$ inversion, as in Case 2 above, i.e. $\pi_1^{**}:A_1(i^{**},j^{**},k^{**}) \to B_1$, $\pi_2^{**}:A_2(i^{**},j^{**},k^{**}) \to B_3$, $\pi_3^{**}:A_3(i^{**},j^{**},k^{**}) \to B_2$.
\item The new triples obtained by applying a horizontal swap to $(\pi_1,\pi_2)$ in $\pi$ and a diagonal swap to $(\pi^{**}_2,\pi_3^{**})$ in $\pi^{**}$ are equal.
%, i.e. $(a_h(\pi_1,\pi_2),\pi_3)=(\pi_1, a_d(\pi_2^{**},\pi_3^{**}))$.
\end{enumerate}
We assume there are such triples, and recall that $i,j,k$ is the sequence corresponding to $\pi$; to shorten notation, let $A_1,A_2,A_3$ be the starting points of the paths determined by $\pi$ and let $A_1^*,A_2^*$ be the initial points of the paths $\pi_1^*:A_1^* \to B_1$ and $\pi_2^*:A_2^* \to B_2$ obtained from the horizontal swap of $(\pi_1,\pi_2)$. Let also $l_1=D_h(\pi_1,\pi_2)$, $l_2=D_d(\pi_2^*,\pi_3)$ and let $(i^*,j^*,k^*)$ be the sequence determining $A_1^*,A_2^*$ and $A_3^*=A_3$. Refer to Figure 3 for the configuration of $A_1,A_2,A_3$. The next lemma shows the relations between $i,i^*,i^{**}$ and so on, needed later.
\begin{lemma}(a) $i^*=i$ and $i^{**}=i+l_2$.\\
(b) $k^*=k+l_1$ and $k^{**}=k^*$.\\
(c) $j^*=j-l_1$ and $j^{**}=j^*$.
\end{lemma}
\begin{proof} This follows from the equations \ref{eqcoordc1} and \ref{eqcoordc2} which record the transformations of $i,j,k$ after a horizontal or diagonal swap.\end{proof}

Since $\pi$ creates an $(12)$ inversion, it must be that $A_2$ is strictly S and strictly E of $A_1$, i.e. \begin{equation}\label{E:inv} x_{A_2} > x_{A_1} \textrm{ and } y_{A_2} < y_{A_1} \end{equation} (use Lemma \ref{lemmaconfig} and the fact that if $x_{A_2} \leqslant x_{A_1}$ then $\pi_1$ and $\pi_2$ must intersect). Moreover, by the definition of $l_1$, \[ x_{A_2}-x_{A_1} \geqslant l_1 \Leftrightarrow (2-i+j)-(k+3) \geqslant l_1 \Leftrightarrow j-i-l_1\geqslant k+1.\]  In particular, \begin{equation}\label{jeq}j-i \geqslant l_1+1.\end{equation} 
\begin{lemma}\label{A} $A_3$ is strictly S and strictly W of $A_2^*$, i.e. $x_{A_2^*}>x_{A_3}$ and $y_{A_2^*}>y_{A_3}$.\end{lemma}
\begin{proof} We have
 \[y_{A_2^*}-y_{A_3}=(\lambda_2 - \lambda_3)+(k-i)+(j-i)+1 \geqslant 2 >0, \] where `$\geqslant$' follows from $\lambda_2 \geqslant \lambda_3$, $k\geqslant i$ and equation (\ref{jeq}). As for $x_{A_2^*}>x_{A_3}$ this happens since $x_{A_2^*}\geqslant x_{A_1}>x_{A_3}$; the first inequality holds because performing a horizontal swap to paths $\pi_1, \pi_2$ creating a $(12)$ inversion implies $A_2^* \in \pi_1$, therefore $A_2^*$ is weakly East of $A_1$; for the second inequality use Lemma \ref{lemmaconfig}.\end{proof}
 %because $\pi$ created an $(12)$ inversion.\end{proof}
 
This lemma, together with the definition of $l_2$, implies that \[ l_2 \geqslant \min \{x_{A_2^*}-x_{A_3},y_{A_2^*}-y_{A_3}\}.\] The triple $(i^{**},j^{**},k^{**})$ must satisfy the constraints (\ref{*}), so in particular
\[i^{**}\leqslant k^{**} \Leftrightarrow i+l_2 \leq k+l_1.\] Then the theorem follows from the following lemma, which contradicts the existence of such a triple, and therefore of $\pi^{**}$.
\begin{lemma} $i +\min \{x_{A_2^*}-x_{A_3},y_{A_2^*}-y_{A_3}\} > k+l_1$.\end{lemma}
\begin{proof} We first show that $i+y_{A_2^*}-y_{A_3} >k+l_1$. This is equivalent to \[ i+(\lambda_2-\lambda_3)+k-i+1+j-i>k+l_1 \Leftrightarrow (\lambda_2-\lambda_3)+1+j-i-l_1>0\] and the last expression is true by equation (\ref{jeq}) above. Similarly, taking into account that $x_{A_2^*}=2-i+j-l_1$, we have \[i+x_{A_2^*}-x_{A_3} >k+l_1 \Leftrightarrow 1+j-l_1+j-i-l_1>0\] which is true again by the equation (\ref{jeq}) above. This finishes the proof of the lemma and of the theorem. \end{proof}

The proof of the Theorem suggests a positive formula to compute $c(\lambda,\mu)$: for a fixed triple $(i,j,k)$, the triples of paths $\pi=(\pi_1,\pi_2,\pi_3)$ which contribute to $c(\lambda,\mu)$ must satisfy the following:\\
\noindent (P1) $\pi$ corresponds to identity permutation, i.e. $\pi_r:A_r \to B_r$, for $1 \leqslant r \leqslant 3$.\\
\noindent (P2) Let $l_h$ be the horizontal distance between $\pi_1$ and $\pi_2$ and let $l_d$ be the diagonal distance between $\pi_2$ and $\pi_3$. Then either $l_h=\infty$,  or, if $l_h < \infty$, neither of triples \[ (i,j+l_h,k-l_h) \textrm{ or } (i+l_d,j,k) \] is triangular, i.e. neither of them satisfies the conditions from (\ref{*}). For the first triple, this can happen, for example, if $j+l_h > \lambda_3$ or if $k < l_h$.

We call the triples $\pi=(\pi_1,\pi_2,\pi_3)$ satisfying (P1) and (P2) {\bf balanced}. The properties (P1) and (P2) mean that one cannot perform a horizontal transformation to $\pi_1$ and $\pi_2$, or a diagonal transformation to $\pi_2$ and $\pi_3$, and obtain a triple of paths with initial points $A_1,A_2,A_3$ coming from a triangular sequence in $S(\lambda)$. Informally, $\pi_2$ is `far enough' from either $\pi_1$ and $\pi_3$, so one cannot do either transformation. In sum:

\begin{cor}\label{posd3} $c(\lambda,\mu)$ is equal to \[ \sum_{s \in S(\lambda)} \# \mathcal{P}_{bal}\bigl((A_1(s),A_2(s),A_3(s)) \to (B_1,B_2,B_3)\bigr) \] where $\mathcal{P}_{bal}$ indicates that only the balanced triples from $\mathcal{P}\bigl((A_1(s),A_2(s),A_3(s)) \to (B_1,B_2,B_3)\bigr)$ are considered. \end{cor}

 \end{document}